\documentclass[10pt]{IEEEtran}
\usepackage{epsfig,color,amsmath,cite}

\IEEEoverridecommandlockouts
\usepackage{wrapfig}
\usepackage{setspace}  
\usepackage{xcolor}
\usepackage{epstopdf}
\usepackage{amssymb}
\usepackage{url}
\usepackage{mathtools}
\usepackage{enumitem}
\usepackage{multirow}
\usepackage{makecell}
\usepackage{amsmath}
\usepackage{stackrel}
\usepackage{booktabs}
\usepackage[flushleft]{threeparttable}
\usepackage{makecell}

\usepackage{algorithm,algorithmic}
\setlist[itemize]{leftmargin=*}

\makeatother

\newcommand{\m}{\boldsymbol}
\allowdisplaybreaks[4]
\usepackage[colorlinks = true,
linkcolor = blue,
urlcolor  = blue,
citecolor = blue,
anchorcolor = blue]{hyperref}

\usepackage{amsthm}
\theoremstyle{definition}
\newtheorem{definition}{Definition}

\newtheorem{property}{Property}
\usepackage[usestackEOL]{stackengine}
\stackMath
\usepackage[framemethod=TikZ]{mdframed}
\usepackage{bm}
\mdfdefinestyle{MyFrame}{%
	linecolor=black,
	outerlinewidth=1.5pt,
	roundcorner=1.5pt,
	innerrightmargin=5pt,
	innerleftmargin=5pt,}
\begin{document}
	\title{\LARGE \bf Geometric Programming-Based Control for Nonlinear, DAE-Constrained Water Distribution Networks}
	\author{Shen Wan$\text{g}^*$, Ahmad F. Tah$\text{a}^*$, Nikolaos Gatsi$\text{s}^{*}$, and Marcio Giacomon$\text{i}^{\dagger}$ \thanks{
			$^*$Department of Electrical and Computer Engineering, $^{\dagger}$Department of Civil and Environmental Engineering, The University of Texas at San Antonio, TX 78249. Emails: mvy292@my.utsa.edu, \{ahmad.taha, nikolaos.gatsis, marcio.giacomoni\}@utsa.edu. This material is based upon work supported by the National Science Foundation under Grant CMMI-DCSD-1728629. }}
	\maketitle
\begin{abstract}
		Control of water distribution networks (WDNs) can be represented by an optimization problem with hydraulic models describing the nonlinear relationship between head loss, water flow, and demand. The problem is difficult to solve due to the non-convexity in the equations governing water flow. Previous methods used to obtain WDN controls (i.e., operational schedules for pumps and valves) have adopted simplified hydraulic models. One common assumption found in the literature is the modification of WDN topology to exclude loops and assume a known water flow direction. 
		In this paper, we present a new geometric programming-based model predictive control approach, designed to solve the water flow equations and obtain WDN controls. The paper considers the nonlinear difference algebraic equation (DAE) form of the WDN dynamics, and the GP approach amounts to solving a series of convex optimization problems and requires neither the knowledge of water flow direction nor does it restrict the water network topology. A case study is presented to illustrate the performance of the proposed method. 
	\end{abstract}

	\section{Introduction and Paper Contributions}~\label{sec:literature}
	Water distribution networks (WDNs) are crucial infrastructures in urban areas. With the expansion of cities, the complexity of WDNs poses challenges for utilities taking into account multiple---potentially conflicting---objectives such as minimizing economic costs, guaranteeing the stability of the network, and maintaining safe water levels in tanks. The very basic decision-making problem involved in the majority of WDN operation problems involves solving for the \textit{water flow} and \textit{head} in pipes given water demand forecasts.
	The hydraulic models of head loss and water flow across pipes, valves, and pumps are accurate, yet highly nonlinear and complex.
	This subsequently hinders optimally solving management/operation problems incorporating the WFP.
	
	The literature of solving the nonconvex WFP is indeed rich and elaborate~\cite{cross1936analysis,wood1972hydraulic,arora1976flows,Hafez-FixedPointWDSA}; a discussion on the merits of these methods is outside the scope of this paper. 
	Most of these methods are iterative algorithms developed to solve a set of linear and nonlinear equations to obtain the physical status  (flows and heads) of WDNs. In their recent paper~\cite{singh2019flow}, the authors investigate the uniqueness of the WFP's solution for generic WDNs.
	
	Several methods have been developed to solve operational and  pump/valve scheduling problems incorporating the WFP in WDNs, and have been recently surveyed in~\cite{mala2017lost} in great detail.
	Model predictive control formulations have been reported in~\cite{wang2017non,sankar2015optimal}. Specifically, the authors in~\cite{wang2016stochastic} present a stochastic MPC formulation to handle uncertainty in WDNs and apply the proposed MPC to the Barcelona drinking water network with real demands.  The authors in~\cite{wang2017non} address a nonlinear economic MPC strategy to minimize the costs associated with pumping and water treatment. A nonlinear MPC controller is designed in~\cite{sankar2015optimal} to meet consumer demands at desired pressures.
	
	The nonlinearity present in WDNs forms a set of nonlinear difference algebraic equations (DAE) which are difficult to handle when solving operation problems. Some recent methods to deal with the nonlinearity are: linearizing the objectives and constraints~\cite{sun2016combining}, constructing relaxations for the nonlinear relationships~\cite{humpola2013unified,singh2018optimal},  keeping the nonlinearity and formulating the problem as a nonconvex program~\cite{xie2015nonlinear}, applying convex approximations to convert the nonconvex problem into a convex one~\cite{sela2015control,fooladivanda2017energy}.
	
	The two works closest to our paper are~\cite{zamzam2018optimal} and~\cite{sela2015control}. The authors in~\cite{zamzam2018optimal} model WDNs as a directed graph, assume the directions of water flow in pipes do not change, and choose the Darcy-Weisbach equation to model head loss in pipes. In~\cite{sela2015control}, the authors convert the nonconvex head loss equations into convex models using geometric programming (GP) approximations, and hence a globally optimal solution is guaranteed. An important contribution of~\cite{sela2015control} is that the proposed GP method is non-iterative. However, it only works under the two assumptions of a tree network topology and known/unchanged flow direction. 
	
	The objective of this paper is to investigate convex optimization-based methods to solve the nonconvex {WFP}, and to subsequently solve an MPC formulation to manage WDN controllers without restricting the WDN graph topology or assuming knowledge of the water flow direction. The proposed approach amounts to solving a series of convex optimization problems, namely, geometric programs, and is embedded within the MPC. We specifically show that the nonconvex, DAE-constrained optimal control problem of WDNs can be approximated by a convex one and solved efficiently. The organization of this paper is given next. Section~\ref{sec:Model} discusses control systems-oriented modeling and MPC formulation of WDNs. Section~\ref{sec:GPmodeling} presents the introduction to GP, and the conversion of the nonconvex hydraulic models to their corresponding convex, GP forms follows. Section~\ref{sec:test} concludes the paper with numerical tests.
	\begin{table*}
		\vspace{-0.3cm}
		\fontsize{7}{7}\selectfont
		\centering
		\renewcommand{\arraystretch}{1.2}
		\caption{WDNs model and its Difference Algebraic Equation (DAE) and Geometric Programming (GP) forms. }
				\vspace{-1em}
		\begin{tabular}{c|c|c|c|c}
			\hline
			& {\textit{Original Hydraulic Model}} & {\textit{DAEs}} \hspace{-2em}& \textit{GP Form}\hspace{-2em} & \textit{Abstracted GP} \\ \hline
			{\hspace{-1em}\textit{Tanks}} \hspace{-2em}
			&
			\parbox{7cm}{
				\vspace{-0.5em}
				\begin{align}~\label{equ:tankhead}
					\hspace{-14pt}h_{i}^{\mathrm{TK}}(k+1)\hspace{-2pt} =\hspace{-2pt} h_{i}^{\mathrm{TK}}(k) \hspace{-2pt}+\hspace{-2pt} \frac{\Delta t}{A_i^{\mathrm{TK}}}\hspace{-3pt}\left(\hspace{-1pt}\sum_{j = 1}^{|\mathcal{N}_i^\mathrm{in}|}\hspace{-3pt}q_{ji}(k)\hspace{-2pt}-\hspace{-2pt}\sum_{j = 1}^{|\mathcal{N}_i^\mathrm{out}|} \hspace{-3pt}q_{ij}(k)\hspace{-3pt}\right)\hspace{-3pt} 
				\end{align}
				\vspace{-0.5em}
			}
			&~\eqref{equ:tankhead-abcstracted} 
			&
			\parbox{6.1cm}{
				\vspace{-0.5em}
				\begin{align}~\label{equ:tankheadNew-exp}
					\hspace{-13pt}{{\hat{h}_i}(k)} \hspace{-3pt}\left(\hspace{-3pt} \prod_{j = 1}^{|\mathcal{N}_i^\mathrm{in}|}\hspace{-3pt}{\hat{q}_{ji}}(k)\hspace{-6pt} \prod_{j = 1}^{|\mathcal{N}_i^\mathrm{out}|}\hspace{-3pt}{\hat{q}_{ij}}^{-1}(k)\hspace{-3pt}\right)\hspace{-4pt}^{\frac{\Delta t}{A_i^{\mathrm{TK}}}}\hspace{-1pt} {{\hat{h}_i}^{-1}\hspace{-2pt}(k\hspace{-2pt}+\hspace{-2pt}1)}\hspace{-2pt}=\hspace{-2pt} 1
				\end{align}
				\vspace{-0.5em}
			}
			&\eqref{equ:tankhead-gp-abcstracted}
			\\
			\hline
			{\textit{\hspace{-5pt}\makecell{Junction \\ nodes}\hspace{-5pt}}} 
			&
			
			\parbox{7cm}{
				\vspace{-0.5em}
				\begin{align}~\label{equ:nodes}
				\textstyle	\sum_{j = 1}^{|\mathcal{N}_i^\mathrm{in}|} q_{ji}(k) - \textstyle \sum_{j = 1}^{|\mathcal{N}_i^\mathrm{out}|} q_{ij}(k) = d_i(k)
				\end{align}
				\vspace{-0.5em}
			}
			& ~\eqref{equ:nodes-abcstracted}
			&
			\parbox{6cm}{
				\vspace{-0.5em}
				\begin{align}~\label{equ:nodes-exp}
				\textstyle	\prod_{j = 1}^{|\mathcal{N}_i^\mathrm{in}|} \hat{q}_{ji}(k)  \textstyle \prod_{j = 1}^{|\mathcal{N}_i^\mathrm{out}|} {\hat{q}_{ij}}^{-1}(k){{\hat{d}_i}^{-1}(k)} = 1
				\end{align}
				\vspace{-0.5em}
			}
			& ~\eqref{equ:node-gp-abcstracted}
			\\
			\hline

			{\hspace{-2em} \textit{Pipes}\hspace{-2em}}
			&  \parbox{7cm}{
				\vspace{-0.7em}
				\begin{align}~\label{equ:head-flow-pipe}
					h_{ij}^\mathrm{P}(k)  = h_{i}(k) - h_{j}(k) = R {q_{ij}(k)}|q_{ij}(k)|^{\mu-1}
				\end{align}
				\vspace{-0.7em}
			}
			& 
			\multirow{2}{*}{\eqref{equ:PumpPipe-abstract}} 
			&
			\parbox{6cm}{
				\vspace{-0.7em}
				\begin{align}~\label{equ:head-loss-pipe-exp}
					{\hat{h}_{i}(k)} {\hat{h}_{i}^{-1}(k)} [C^{\mathrm{P}}(k)]^{-1} {\hat{q}_{ij}}^{-1}(k) = 1
				\end{align}
				\vspace{-0.7em}
			} 
			& ~\eqref{equ:Pipe-gp-abstract}
			\\ \cline{1-2} \cline{4-5}
			
			{\hspace{-2em} \textit{Pumps}\hspace{-2em}} 
			& \parbox{7cm}{
				\vspace{-0.7em}
				\begin{align} \label{equ:head-flow-pump}
					\hspace{-10pt}h_{ij}^\mathrm{\mathrm{M}}(k) = h_{i}(k) - h_{j}(k) = -{s_{ij}^2(k)}(h_0 - r  (q_{ij} s_{ij}^{-1})^\nu )
				\end{align}
				\vspace{-0.7em}
			} 
			&
			&\parbox{6cm}{
				\vspace{-0.7em}
				\begin{align} \label{equ:head-flow-pump-exp}
					\hspace{-12pt}{\hat{h}_{i}(k)} {\hat{h}_{i}^{-1}(k)}{{\hat{s}_{ij}}(k)}^{-C_1^{\mathrm{M}}(k)}{{\hat{q}_{ij}}(k)}^{-C_2^{\mathrm{M}}(k)}  = 1
				\end{align}
				\vspace{-0.7em}
			}
			& ~\eqref{equ:Pumps-gp-abcstracted}
			\\
			\hline
			
			{\hspace{-3.5em} \textit{Constraints}\hspace{-2em}} 
			
			& \parbox{7cm}{
				\vspace{-1em}
				\begin{subequations} ~\label{equ:constraints}
					\begin{align}
						h_{i}^{\mathrm{min}} &\leq  h_{i}(k) \leq h_{i}^{\mathrm{max}}~\label{equ:tankLimit} \\
						&h_i^{\mathrm{R}}(k) = h_i^{\mathrm{R}}~\label{equ:tank-reservoir}\\
						0 &\leq  s_{ij}(k) \leq 1 ~\label{equ:speedLimit} \\
						q_{ij}^{\mathrm{min}} &\leq  q_{ij}(k) \leq q_{ij}^{\mathrm{max}} ~\label{equ:flowLimit}
					\end{align}
				\end{subequations}
				\vspace{-1em}
			} 
			& \eqref{equ:constr-abcstracted}
			&
			
			\parbox{6cm}{
				\vspace{-0.7em}
				\begin{subequations}~\label{equ:tankLimit-exp}
					\begin{align}
						\hat{h}_{i}^{-1}(k)	\hat{h}_{i}^{\mathrm{min}} \leq  1,{\hat{h}_{i}(k)}  \left( \hat{h}_{i}^{\mathrm{max}}\right) ^{-1} \leq 1 \\
						{\hat{h}_i}^{-1}(k)  {\hat{h}_i}^{\mathrm{R}} = 1 ~\label{equ:tank-reservoir-exp} \\
						{\hat{s}_{ij}}^{-1}(k)	 \leq  1,{\hat{s}_{ij}}(k)  b^{-1} \leq 1 \\
						{\hat{q}_{ij}}^{-1}(k)	{\hat{q}_{ij}}^{\mathrm{min}} \leq  1,{\hat{q}_{ij}}(k)  \left( {\hat{q}_{ij}}^{\mathrm{max}}\right) ^{-1} \leq 1 
					\end{align}
				\end{subequations}
				\vspace{-1em}
			}
			& \eqref{equ:constr-gp-physical}
			\\
			\hline \hline
		\end{tabular}
		\label{tab:models}%
		\vspace{-1.5em}
	\end{table*}
	
	\vspace{-0.25cm}
	\section{Modeling of WDNs}~\label{sec:Model}
	We model WDNs by a directed graph $(\mathcal{V},\mathcal{E})$.  Set $\mathcal{V}$ defines the nodes and is partitioned as $\mathcal{V} = \mathcal{J} \bigcup \mathcal{T} \bigcup \mathcal{R}$ where $\mathcal{J}$, $\mathcal{T}$, and $\mathcal{R}$ stand for the collection of junctions, tanks, and reservoirs. Let $\mathcal{E} \subseteq \mathcal{V} \times \mathcal{V}$ be the set of links, and define the partition $\mathcal{E} = \mathcal{P} \bigcup \mathcal{M} \bigcup \mathcal{W}$, where $\mathcal{P}$, $\mathcal{M}$, and $\mathcal{W}$ stand for the collection of pipes, pumps, and valves. For the $i$-th node, set $\mathcal{N}_i$ collects its neighboring nodes and is partitioned as $\mathcal{N}_i = \mathcal{N}_i^\mathrm{in} \bigcup \mathcal{N}_i^\mathrm{out}$, where $\mathcal{N}_i^\mathrm{in}$ and $\mathcal{N}_i^\mathrm{out}$ stand for the collection of inflow and outflow nodes. It is worth emphasizing that the assignment of direction to each link 
is arbitrary, as the  presented optimizations yield optimal flow direction in pipes.
	 
	\subsection{Conservation of mass and energy}~\label{sec:Model_iass}
	The hydraulic equations describing the flow in  WDNs are derived from the principles of \textit{conservation of mass} and \textit{energy}. 
	According to these basic laws, the equations that model mass and energy conservation for all components (tanks, junctions, pipes, and pumps) in WDNs can be written and provided in explicit and compact matrix-vector forms in the first three columns of Table~\ref{tab:models}. The last two columns of Table~\ref{tab:models} are needed in the ensuing sections of the paper.
	
	\noindent \textbf{Tanks and Junction Nodes} --- The water volume evolution in the $i$-th tank at time $k$ can be expressed by a discrete-time difference equation~\eqref{equ:tank-volume}, while the head created by the tank can be described as~\eqref{equ:tank-head}
	
	\vspace{-1.2em}
	\begin{subequations}
		{\small	\begin{eqnarray}
			V_{i}(k+1) &\hspace{-1em}=\hspace{-1em}& V_{i}(k) + \Delta t \hspace{-3pt}\left(\hspace{-1pt}\textstyle \sum_{j = 1}^{|\mathcal{N}_i^\mathrm{in}|} q_{ji}(k)\hspace{-2pt}-\hspace{-2pt} \textstyle \sum_{j = 1}^{|\mathcal{N}_i^\mathrm{out}|} q_{ij}(k)\hspace{-2pt}\right)\hspace{-3pt}  ~\label{equ:tank-volume} \\
			h_i^{\mathrm{TK}}(k) &\hspace{-1em}=\hspace{-1em}& \frac{V_i(k)}{A_i^{\mathrm{TK}}} + E_i^{\mathrm{TK}}, \; i\in\mathcal{T}, ~\label{equ:tank-head} 
			\end{eqnarray}}
	\end{subequations}
	\vspace{-1em}

	\noindent where $V_{i}$ and $\Delta t$ are the volume and sampling time; $q_{ji}(k),\;i \in \mathcal{J},\;j \in \mathcal{N}_i^\mathrm{in} $ stands for the inflow of $j$-th neighbor, while $q_{ij}(k),\;i \in \mathcal{J},\;j \in \mathcal{N}_i^\mathrm{out} $ stands for the outflow of $j$-th neighbor; $h_i^{\mathrm{TK}}$, $A_i^{\mathrm{TK}}$, and $E_i^{\mathrm{TK}}$ stand for the head, cross-sectional area, and elevation of the $i$-th tank.
	Combining~\eqref{equ:tank-volume} and~\eqref{equ:tank-head}, the head changes from time $k$ to $k+1$ of the $i$-th tank can be written as~\eqref{equ:tankhead} in Table~\ref{tab:models}. 
	Junction nodes are the points where water flow merges or splits. The expression of mass conservation of the $i$-th junction at time $k$ can be written as~\eqref{equ:nodes} in Table~\ref{tab:models}, and $d_i(k)$ stands for the end-user demand that is extracted from node $i$.
	
	\noindent \textbf{Head Loss in Pipes} --- The major head loss of a pipe from node $i$ to $j$ is due to friction and is determined by~\eqref{equ:head-flow-pipe} from Table~\ref{tab:models}, where $h_i$ and $h_j$ stand for the head of the $i$-th and $j$-th node; $h_{ij}^{\mathrm{P}}$ is the head loss of the pipe from $i$ to $j$; $R$ is the resistance coefficient and $\mu$ is the constant flow exponent in the corresponding formula. Table~\ref{tab:headloss} represents the most common formulae used in the literature to model the resistance coefficient $R$.  The approach presented in this paper considers any of the three formulae in Table~\ref{tab:headloss}. 
	The numerical tests use the Hazen-Williams formula~\cite{rossman2000epanet}, which is typically used to analyze urban water supply systems.

	\noindent \textbf{Head Gain in Pumps/Reservoirs} --- A head increase can be generated by a pump between the suction node $i$ and the delivery node $j$. The pump properties decide the relationship function between pump flow and head increase~\cite{rossman2000epanet}. Generally, it can be expressed as~\eqref{equ:head-flow-pump}, where $h_{ij}^{\mathrm{M}}$ is head increase of the pump from $i$ to $j$;  $h_0$ is the shutoff head for the pump; $q_{ij}$ is the flow through a pump; $s_{ij} \in [0,1]$ is the relative speed of the same pump; $r$ and $\nu$ are the pump curve coefficients.   
	We also assume that reservoirs have infinite water supply, the head of the $i$-th reservoir $h_i^{\mathrm{R}}$ is fixed, and this can be viewed as an operational constraint.
	
	\noindent \textbf{WDN Operational Constraints} ---
	The constraints (upper and lower bounds on the heads of junctions, tanks and reservoirs, pump speeds, and flows) are expressed as  equations~\eqref{equ:tankLimit}--\eqref{equ:flowLimit} in Table~\ref{tab:models}. Note that constraint~\eqref{equ:speedLimit} can be replaced with binary (on/off) constraints on the pump schedules if variable speed pumps are not available in the water network. In this paper, we assume that the relative speed of all pumps can be modulated in the interval $[0,1]$, rather that needing integer constraints. 

	\begin{table}[t]
		\fontsize{7}{7}\selectfont
		\caption{Head loss formula$\text{e.}^*$}
		\begin{threeparttable}
			\centering
			\makegapedcells
			\setcellgapes{1pt}
			\begin{tabular}{ c|c|c }
				\hline
				\textit{Formula} & \textit{Resistance Coefficient} ($R$) & \textit{Flow Exponent} ($\mu$) \\ \hline
				Hazen-Williams &  $4.727 L^{\mathrm{P}} C_{\mathrm{HW}}^{-1.852} (D^{\mathrm{P}})^{-4.871}$ & 1.852   \\ \hline
				Darcy-Weisbach &  $0.0252 L^{\mathrm{P}} f(\epsilon,D^{\mathrm{P}},q) (D^{\mathrm{P}})^{-5}$ & 2   \\ \hline
				Chezy-Manning &  $4.66  L^{\mathrm{P}} C_{\mathrm{CM}}^{2} (D^{\mathrm{P}})^{-5.33}$ & 2   \\ \hline
				\hline
			\end{tabular}
			\begin{tablenotes}
				\scriptsize
				\item $^*C_{\mathrm{HW}}$, $\epsilon$, $C_{\mathrm{CM}}$ are  roughness coefficients of Hazen-Williams, Darcy-Weisbach and Chezy-Manning. $D^{\mathrm{P}}\ (\mathrm{ft})$ is the pipe diameter, $L^{\mathrm{P}} \ (\mathrm{ft})$ is the pipe length. $q \ (\mathrm{cfs})$ is the flow rate, $f$ is friction factor (dependent on $\epsilon$, $D^{\mathrm{P}}$, and $q$).
				\vspace{-0.15cm}	
		\end{tablenotes}
			\label{tab:headloss}
					\vspace{-0.15cm}
		\end{threeparttable}
			\vspace{-0.15cm}
	\end{table}
			\vspace{-0.245cm}
	\subsection {Difference algebraic equations (DAE) form of WDNS}
	The WDN model in previous section can be abstract in the form of DAEs as~\eqref{equ:dae-abstract}. Define $\m x$, $\m l$, $\m u$, $\m v$, and $\m s$ to be vectors collecting heads at tanks, heads at junctions, flows through pumps and valves, flows through pipes, and the relative speed of pumps. And the corresponding dimensions are $n_t$,$n_j$,$n_{w}+n_{m}$,$n_v$, and $n_m$. Collecting the mass and energy balance equations of tanks~\eqref{equ:tankhead}, junctions~\eqref{equ:nodes}, pipes~\eqref{equ:head-flow-pipe}, and pumps~\eqref{equ:head-flow-pump}, we write the following DAE model
	\begin{subequations} ~\label{equ:dae-abstract}
		\begin{align}
			\textit{DAE:}\;\;\;	\m x(k + 1) &= \m A \m x(k) + \m B_u \m u(k) + \m B_v \m v(k) ~\label{equ:tankhead-abcstracted} \\
			\hspace{-10pt}	\m 0_{n_{j}} &=\m E_u \m u(k) + \m E_v \m v(k) + \m E_d \m d(k)  ~\label{equ:nodes-abcstracted}\\
			\hspace{-10pt}\m 0_{n_{w}+n_{m}}&= \m E_x \m x(k) + \m E_l \m l(k) + \m \Phi(\m u, \m v, \m s)    ~\label{equ:PumpPipe-abstract}
		\end{align}
	\end{subequations}
	where $\m A$, $\m E$, and $\m B$ are constant matrices that depend on the WDN topology and the aforementioned hydraulics. The function $\m \Phi(\cdot)$
collects the nonlinear components in \eqref{equ:head-flow-pipe} and~\eqref{equ:head-flow-pump}.
	The physical constraints pertaining to~\eqref{equ:dae-abstract} can be rewritten as
	
	\vspace{0.15cm}
	
	\noindent \textit{Constraints:}
	\begin{align}
		&\m x(k) \in [\m x^{\mathrm{min}}(k),\m x^{\mathrm{max}}(k)],\m l(k) \in [\m l^{\mathrm{min}}(k),\m l^{\mathrm{max}}(k)] \notag \\
		&\m u(k) \in [\m u^{\mathrm{min}}(k),\m u^{\mathrm{max}}(k)],\m v(k) \in [\m v^{\mathrm{min}}(k),\m v^{\mathrm{max}}(k)]  \notag  \\
		&\m s(k) \in [\m 0_{n_m},\m 1_{n_m}] ~\label{equ:constr-abcstracted} .
	\end{align}
The next section discusses the WDN control objectives.
 
	\subsection {Control objectives}~\label{sec:CtrlObj}
	Multiple objectives can be applied depending on operational considerations, and can be expressed through
%
	\begin{align}
	&\hspace{-5pt} \Gamma_1(\m x(k))\hspace{-2pt}=\hspace{-2pt} \left\{\Centerstack[l]{ (\m x(k)\hspace{-2pt} -\hspace{-2pt} \m x^{\mathrm{sf}})^{\top}(\m x(k)\hspace{-2pt} -\hspace{-2pt} \m x^{\mathrm{sf}}),\\ 0}\right.&&
	\Centerstack{\hspace{-2pt} \mathrm{if}\ \m x(k)\hspace{-2pt} \leq \hspace{-2pt} \m x^{\mathrm{sf}} \\\hspace{-2pt} \text{Otherwise} } ~\label{equ:SafetyWater} \\
	&\hspace{-5pt} \Gamma_2(\Delta \m u(k)) = \Delta \m u(k)^{\top}\Delta \m u(k) &&   ~\label{equ:SmoothinControl}
	\end{align}
	where $\Gamma_1(\cdot)$ enforces maintaining the safety water storage in each tankdecided by the operator and $\m x^{\mathrm{sf}}$ is a vector collecting the safety head levels of tanks; $\Gamma_2(\cdot)$ enforces the smoothness of control actions through $\Delta \m u(k) = \m u(k) - \m u(k-1)$ which stands for the flow rate changes of controllable components from time $k-1$ to $k$. At initial time $k=t_0$, define a vector that collects all the optimization variables as
	$$ \boldsymbol \xi [t_0] \triangleq \Bigl \lbrace \m x(k+1),\m u(k), \m l(k), \m v(k), \m s(k) \Bigr \rbrace_{k=t_0}^{k=t_0+H_p}, $$
	where $H_p$ is the prediction horizon of the MPC. {Note that the indexing for $\m x(k)$ is different in $\boldsymbol \xi[t_0]$ due to the fact that the initial conditions of the tanks $\m x(t_0)$ is known, unlike other optimization variables such as the flow and the pump controls which we need to solve for from $k=t_0$ through $k=t_0+H_p$.} The multi-objective cost function can be written as
	\begin{equation*}~\label{equ:Multi-objective}
		\Gamma(\boldsymbol \xi [t_0] ) =\textstyle \sum_{i=1}^{2}\Gamma_{i}(\boldsymbol \xi[t_0])= \boldsymbol\xi^{\top}[t_0] \m \Omega \boldsymbol\xi[t_0]+ \m \omega^{\top} \boldsymbol\xi[t_0] + \omega,
	\end{equation*}
	where $\m \Omega, \m \omega$ and $\omega$ are the corresponding weight matrix, vector, and scalars from $\Gamma_{1,2}(\boldsymbol \xi)$. Similar objective functions have been used before in~\cite{ocampo2013application,sun2016combining}.
%

	\subsection {MPC formulation}
	Here, we propose using MPC to solve the WDN operation problem considering the nonlinearities and nonconvexities present in the energy balance equations in the DAE model~\eqref{equ:dae-abstract}. At each time instant $k$, MPC obtains the control input to be implemented, given a prediction of the current disturbance, output, and states, and then finds an optimal control sequence over a prediction horizon $H_p$. 
	The MPC can be written as
		\begin{align}
	\min_{\boldsymbol \xi[t_0]} \;\;& \Gamma\left(\boldsymbol \xi [t_0]  \middle|\: \m x(t_0), \left\lbrace \m d(k) \right\rbrace_{k=t_0}^{k=t_0+H_p} \right) \notag \\
			\mathrm{s.t.}\;\;& \textit{DAE}~\eqref{equ:dae-abstract},\; \textit{Constraints}~\eqref{equ:constr-abcstracted}  \label{equ:NMPC}
		\end{align}
	The optimization problem~\eqref{equ:NMPC} is nonlinear and nonconvex due to the head loss models of pipes and pumps. WDN-MPC requires the knowledge or prediction of the nodal water demand for the entire prediction horizon as well as the initial levels of tanks $\m x(t_0)$, and $\left\lbrace \m d(k) \right\rbrace_{k=t_0}^{k=t_0+H_p}.$
	
	%

	\vspace{-0.25cm}
	
	\section{GP Modeling of WDNs}~\label{sec:GPmodeling}
	A basic introduction to GP is provided first, followed by the conversion of the nonconvex hydraulic models in~\eqref{equ:NMPC} to their corresponding convex, GP form. 	A geometric program (GP) is a type of optimization problem with objective and constraint functions that are monomials and posynomials~\cite{boyd2007tutorial}.
	Here, we propose a GP model without assuming a known flow direction by mapping the optimization variable $\m \xi[t_0]$ in~\eqref{equ:NMPC} into its exponential form. The conversion helps to map all of the non-positive values into positive ones. 
	Specifically, we convert the head and demand at the $i$-th node $h_i$ and $d_i$, the flow $q_{ij}$, and relative speed $s_{ij}$ into positive values ${\hat{h}_i}$, ${\hat{d}_i}$, ${\hat{q}_{ij}}$, and ${\hat{s}_{ij}}$ through exponential functions, as follows
	\begin{equation}~\label{equ:NLPGP}
	{\hat{h}_i} \triangleq {b}^{h_i}, \; {\hat{d}_i} \triangleq {b}^{d_i},\;
	{\hat{q}_{ij}} \triangleq {b}^{q_{ij}} , \; {\hat{s}_{ij}} \triangleq {b}^{s_{ij}},  
	\end{equation}
	where $b=1+\delta$ is a constant base and $\delta$ is a small positive number. After conversion of variables, the ${\hat{h}_i}$, ${\hat{d}_i}$, ${\hat{q}_{ij}}$, and ${\hat{s}_{ij}}$ are positive values which can then be used to transform the nonconvex problem~\eqref{equ:NMPC} into a GP.
	Converting the junction and tank physical models as well as constraints follows from the above exponential mapping~\eqref{equ:NLPGP}, while converting the pipe and pump models into GP form is more complicated. The last two columns of Table~\ref{tab:models} produce detailed as well as abstract versions of the conversions of all physical models. 
	\vspace{-1em}
	\subsection{ Conversion of mass and energy balance equations}~\label{sec:modelconversion}
	For {the models of tanks and  junctions}, the conversion process is straightforward. After exponentiating both sides of~\eqref{equ:tankhead} and~\eqref{equ:nodes}, variables $q_{ij}$, $h_i$, and $d_i$ are changed  into ${\hat{q}_{ij}}$, ${\hat{h}_i}$, and ${\hat{d}_i}$, while constraints~\eqref{equ:tankhead} and~\eqref{equ:nodes} are converted to monomial equality constraints~\eqref{equ:tankheadNew-exp} and~\eqref{equ:nodes-exp} in Table~\ref{tab:models}.
	
	In order to clearly show the derivation for pipes, the time $k$ is ignored at first. At time $k$, let ${\hat{h}_{ij}^{\mathrm{P}}}$ be the GP form of head loss of a pipe, which is obtained  by exponentiating both sides of~\eqref{equ:head-flow-pipe} as follows
	\begin{align*}
	{\hat{h}_{i}} {\hat{h}_{i}^{-1}} = {\hat{h}_{ij}^{\mathrm{P}}}  &= {\large b^{\left(q_{ij} R {|q_{ij}|}^{\mu-1} - q_{ij} \right)}} \cdot b^{q_{ij}} = C^{\mathrm{P}}(q_{ij})\cdot{\hat{q}_{ij}},
	\end{align*}
	where  $C^{\mathrm{P}}(q_{ij})=b^{q_{ij} \left(R {|q_{ij}|}^{\mu-1} - 1\right)}$ is a function of $q_{ij}$.   
	
	Note that variable ${q_{ij}(k)}$ is unknown at each time $k$. The premise is to solve a series of convex optimization problems to find the final value for each time $k$. Hence, we can make an initial guess denoted by $\langle{q_{ij}}\rangle_0$ and  $\langle{C^{\mathrm{P}}}\rangle_0$ for iteration $0$. For the $n$-th iteration, the corresponding values are denoted by $\langle{q_{ij}}\rangle_n$ and $\langle{C^{\mathrm{P}}}\rangle_n$. If the flow rates are close to each other between the $(n-1)$-th and $n$-th iteration, we can approximate $\langle{C^{\mathrm{P}}}\rangle_n$ using  $\langle{C^{\mathrm{P}}}\rangle_{n-1}$, that is
$		\langle{C^{\mathrm{P}}}\rangle_n \approx \langle{C^{\mathrm{P}}}\rangle_{n-1}.$
	Thus, for each iteration $n$,  $$\langle{C^{\mathrm{P}}(k)}\rangle_n = b^{\langle{q_{ij}(k)}\rangle_{n-1} \left(R {|\langle{q_{ij}(k)}\rangle_{n-1}|}^{\mu-1} - 1\right)}$$ can be approximated by a constant given the flow value $\langle{q_{ij}(k)}\rangle_{n-1}$ from the previous iteration. With this approximation, the head loss constraint for each pipe can be written as a \emph{monomial equality constraint}
	${\hat{h}_{i}(k)} {\hat{h}_{i}^{-1}(k)} = C^{\mathrm{P}}(k) {\hat{q}_{ij}}(k)$
	which is equivalently expressed as~\eqref{equ:head-loss-pipe-exp}. 
	If we solve flow $q_{ij}(k)$, this value can be as an initialization for the next iteration, implying that $\langle{C^{\mathrm{P}}(k+1)}\rangle_0 = \langle{C^{\mathrm{P}}(k)}\rangle_n$ which can accelerate the convergence of the successive convex approximation.
	
	Similarly for pumps, the new variables ${\hat{q}_{ij}}(k) = b^{q_{ij}(k)}$ and  ${\hat{s}_{ij}}(k) = b^{s_{ij}(k)}$ for $(i,j) \in \mathcal{M}$ are introduced.  Let ${\hat{h}_{ij}^\mathrm{M}} $ be the GP form of head increase of a pump:
	\begin{align*}
		{\hat{h}_{i}} {\hat{h}_{i}^{-1}}\hspace{-2pt} = \hspace{-2pt}{\hat{h}_{ij}^\mathrm{M}}
		\hspace{-2pt}=\hspace{-2pt}(b^{s_{ij}})^{-{s_{ij}} h_0}\cdot(b^{q_{ij}})^{ r q_{ij}^{\nu-1} s_{ij}^{2-\nu}}
		\hspace{-3pt}=\hspace{-2pt}({\hat{s}_{ij}})^{C_1^{\mathrm{M}}} ({\hat{q}_{ij}})^{C_2^{\mathrm{M}}},
	\end{align*}
	where $C_1^{\mathrm{M}} = -{s_{ij}} h_0$ and $C_2^{\mathrm{M}} = r q_{ij}^{\nu-1} s_{ij}^{2-\nu}$. Parameters $C_1^{\mathrm{M}}(k)$ and $C_2^{\mathrm{M}}(k)$ follow a similar iterative process as $C^{\mathrm{P}}(k)$. That is, they are treated at iteration $n$ as constants based on the flow and relative speed values at iteration $n-1$. Hence, the approximating equation for the pump head increase becomes the monomial equality constraint~\eqref{equ:head-flow-pump-exp}, where $\nu$ is a constant parameter determined by the pump curve. 
	
	Therefore, starting with an initial guess for the flow rates and relative speeds, the constraints in the DAE model are approximated at every iteration by constraints abiding by the GP form, as listed in Table~\ref{tab:models}. 
	This process continues until a termination criterion is met. The details are further discussed in Algorithm~\ref{alg:seq}, after the presentation of the  abstract GP form and the conversion of the control objectives in the next section.

	\vspace{-0.35cm}
	
	\subsection {Abstract GP model}
	To express the GP-based form of WDN-MPC in a compact form, we define some mathematical operators.
	\begin{definition}
		For matrices $\m X$ and $\m B \in \mathbb{R}^{m\times n}$, the \textit{element-wise exponential} operation on $\m X$ with base $\m B$ is a matrix of the same dimension with elements given by
		$$\hat{\m X} = \m B^ {\m X} = \begin{bmatrix}
		b^{x_{11}}_{11} & \cdots & b^{x_{1n}}_{1n} \\
		\vdots & \ddots & \vdots \\
		b^{x_{m1}}_{m1} & \cdots & b^{x_{mn}}_{mn} 
		\end{bmatrix} = \begin{bmatrix}
		{\hat{x}_{11}} & \cdots & {\hat{x}_{1n}} \\
		\vdots & \ddots & \vdots \\
		{\hat{x}_{m1}} & \cdots & {\hat{x}_{mn}} 
		\end{bmatrix}.$$
	\end{definition}
	%
	\begin{definition}
		For matrices $\m Y \in \mathbb{R}^{n\times m}$ and matrix $\m X\in \mathbb{R}^{m\times p}$, the \textit{element-wise exponential matrix product}  $\m C = \m X {\tiny \star} \m Y$ has elements given by $c_{ij} = \prod_{k=1}^{m} (\hat{x}_{kj})^{y_{ik}}$ for $i = 1,\ldots,n$ and $j =1,\ldots,p$, where $\hat{x}_{kj} = b^{{x}_{kj}}$.
	\end{definition}
	\begin{property}~\label{prp1}
		For matrices $\m Y$ with size  $n\times m$ and $\m X$ with size $m\times p$, let $ \hat{\m X} = b^ {\m X}$, where $b$ is base. The following  holds: 
		$$b^{\m Y \m X} = \hat{\m X} {\tiny \star} \m Y.$$
	\end{property}
	
	With the above definitions, we now derive the GP-based DAE model and constraints of WDN. Performing an {element-wise exponential} operation on both sides of~\eqref{equ:dae-abstract} yields
	\begin{align*}
		b^{ \m x(k + 1)} &= b ^{ \m A \m x(k) + \m B_u \m u(k) + \m B_v \m v(k) }\\
		\m 1 &= b ^ {\m E_u \m u(k) + \m E_v \m v(k) + \m E_d \m d(k)} \\
		\m 1 &= b ^ {\m E_x \m x(k) + \m E_l \m l(k) + \m \Phi(\m u(k), \m v(k), \m s(k)) }.
	\end{align*}
	Denote $\m \hat{\m x}(k) = b ^{\m x(k)}$ and similarly $\m l(k)$, $\m u(k)$, $\m v(k)$, and $\m s(k)$  are converted into $\m \hat{\m l}(k)$,$\m \hat{\m u}(k)$,$\m \hat{\m v}(k)$ and $\m \hat{\m s}(k)$. 
	The models of junctions and tanks can be written as monomials~\eqref{equ:tankhead-gp-abcstracted} and~\eqref{equ:node-gp-abcstracted}  directly according to Property~\ref{prp1}.
	
	For a pipe from node $i$ to $j$, according to~\eqref{equ:head-loss-pipe-exp}, the exponential of nonlinear function is $C^{\mathrm{P}}(k){\hat{q}_{ij}}(k)$.  The head loss constraints can be compactly written for all pipes using the element-wise product $\m F_v(k) {\circ} \m \hat{\m v}(k)$, where $\m F_v(k)$ is a $n_p \times 1$ column vector collecting the $C^{\mathrm{P}}(k)$ parameters of pipes.  
	
	Similarly, for all pumps, define $\m F_s(k)$ and $\m F_u(k)$ as $n_m \times 1$ column vectors respectively collecting all of parameters $C_1^{\mathrm{M}}(k)$ and $C_2^{\mathrm{M}}(k)$ of pumps. In summary,~\eqref{equ:Pipe-gp-abstract} and~\eqref{equ:Pumps-gp-abcstracted} are the abstract GP form of pipe and pump models. The overall DAE-GP model is given as follows
	\vspace{0.2cm}
	
	\noindent	\textit{DAE-GP:}
	\begin{subequations}~\label{equ:gp-abcstracted}
		\begin{align}
			&\hspace{-8pt} \hat{\m x}(k+1) =  [\m A {\tiny \star}\m \hat{\m x}(k)]{\circ}[\m B_u {\tiny \star}\m \hat{\m u}(k)]{\circ}[\m B_v {\tiny \star} \m \hat{\m v}(k)] ~\label{equ:tankhead-gp-abcstracted} \\
			&\hspace{16pt} \m 1_{n_{j}} = [\m E_u{\tiny \star}\m \hat{\m u}(k) ]{\circ}[\m E_v {\tiny \star} \m \hat{\m v}(k) ]{\circ} [\m E_d {\tiny \star} \m \hat{\m d}(k)] ~\label{equ:node-gp-abcstracted} \\
			&\hspace{-8pt} [\m E_x {\tiny \star}\m \hat{\m x}(k)]{\circ}[\m E_l {\tiny \star}\m \hat{\m l}(k)] = \m F_v(k) {\circ} \m \hat{\m v}(k) ~\label{equ:Pipe-gp-abstract}\\
			& \hspace{-8pt} [\m E_x  {\tiny \star} \m \hat{\m x}(k)]{\circ}[\m E_l {\tiny \star} \m \hat{\m l}(k)] = [\m \hat{\m s}(k) ^{ \m F_s(k)}] {\circ} [\m \hat{\m u}(k) ^ {\m F_u(k)}].~\label{equ:Pumps-gp-abcstracted} 
		\end{align}
	\end{subequations}
In addition, the  physical constraints~\eqref{equ:constr-abcstracted} can be rewritten as 
	
	\vspace{0.2cm}
	
	\noindent			\textit{Constraints-GP:}
	\begin{align}
		&\hat{\m x}(k) \in [\hat{\m x}^{\mathrm{min}}(k),\hat{\m x}^{\mathrm{max}}(k)],\hat{\m l}(k) \in [\hat{\m l}^{\mathrm{min}}(k),\hat{\m l}^{\mathrm{max}}(k)]  \notag \\
		&\hat{\m s}(k) \in [\m 1_{n_m},\m b_{n_m}] ~\label{equ:constr-gp-physical} \\
		&\hat{\m u}(k) \in [\hat{\m u}^{\mathrm{min}}(k),\hat{\m u}^{\mathrm{max}}(k)],\hat{\m v}(k) \in [\hat{\m v}^{\mathrm{min}}(k),\hat{\m v}^{\mathrm{max}}(k)].  \notag 
	\end{align}
	
		\vspace{-0.351cm}
	\subsection{Conversion of control objectives}~\label{sec:objconversion}
In this section, we covert the control objectives in the nonconvex problem \eqref{equ:NMPC} to their convex, GP-based form.	In~\eqref{equ:SafetyWater}, notice that $\m x$ is a vector collecting the head $h_i$ at tanks. The objective
	$(\m x(k) - \m x^{\mathrm{sf}})^{\top}(\m x(k) - \m x^{\mathrm{sf}})$ 
	encourages $\m x(k)$ to be close to the constant $\m x^{\mathrm{sf}}$. 
	Hence, we introduce a new auxiliary variable $\hat{\m z}(k) \triangleq b^{\m x^{\mathrm{sf}} - \m x(k)}$ which will be pushed to be close to $\m 1$. Using the epigraph form, the original objective function $\Gamma_1(\m x(k))$ in Section~\ref{sec:CtrlObj} is replaced by
	$\hat{\Gamma}_1(\hat{\m z}(k)) = \prod_{i = 1}^{n_t}  \hat{z}_i(k)$
	and $2n_t$ constraints are added as
	\begin{subequations}~\label{eq:constraints-gp-obj1}
		\begin{align}
	\textit{Safety-GP:} \;\;			& \hat{\m z}(k) \geq \m 1_{n_t}   \\
			& \hat{z}_i(k) = \hat{x}_i^{\mathrm{sf}} \hat{x}_i^{-1}(k)\;\;, i \in [1,n_t]
		\end{align}
	\end{subequations}
	where $\hat{\m x}^{\mathrm{sf}}$ and $\hat{\m x}(k)$ are the GP form of ${\m x}^{\mathrm{sf}}$ and ${\m x}(k)$. If the water level of the $i$-th tank is above the safe level, then variable $\hat{z}_i(k)$ already satisfies the constraint, which implies no objective function is applied. Otherwise, the objective is $\hat{z}_i(k)$, and the corresponding constraint is $\hat{z}_i(k) =\hat{x}_i^{\mathrm{sf}} \hat{x}_i^{-1}(k)$.
	Moving to the second part of the objective function~\eqref{equ:SmoothinControl}, $\Delta \m u(k) = \m u(k) - \m u(k-1)$ is a vector collecting the flow changes of controllable flow $\m u(k)$.  We introduce a new auxiliary variable $\hat{\m p}(k) \triangleq b^{\m u(k)- \m u(k-1)}$ and perform an {element-wise exponential} operation on both sides of~\eqref{equ:SmoothinControl} yielding
	\begin{align*}
		b^{[\m u(k) - \m u(k-1)]^{\top}\Delta \m u(k)} &= [\hat{\m u}(k) \circ \hat{\m u}^{-1}(k-1)]^{\Delta \m u(k)}\\ &= (\hat{p}(k))^{\Delta \m u(k)}.
	\end{align*}
	Using the epigraph form, the original objective function $\Gamma_2(\Delta \m u(k)) $ can be expressed as a new objective
$
		\hat{\Gamma}_2(\hat{\m p}(k)) =\textstyle  \prod_{i = 1}^{n_u}  (\hat{p}_i(k))^{\Delta u_i(k)}$
	and $n_u + 1$ constraints given as
	
	\vspace{0.2cm}
	\noindent 	\textit{Smoothness-GP:}
		\begin{align}~\label{eq:constraints-gp-obj2}
	\hspace{-0.35cm}	 \hat{\m p}_i(k) = \hat{\m u}_i(k) \hat{\m u}^{-1}_i(k-1), \forall \,  i \in \{1,n_u\},\;
		 \hat{\Gamma}_2(\hat{\m p}(k)) \geq 1.
		\end{align}
	Given the above derivations, the final GP form of multi-objective cost function can be rewritten as
	\begin{equation}~\label{eq:gp-obj}
	\hat{\Gamma}(\hat{\m z}(k),\hat{\m p}(k)) =\hat{\Gamma}_1(\hat{\m z}(k)) + \hat{\Gamma}_2(\hat{\m p}(k)),
	\end{equation}
	where $\hat{\Gamma}(\hat{\m z}(k),\hat{\m p}(k))$ is a posynomial function.
The convex GP-based MPC can now be expressed as
		\begin{align}
			\min_{\substack{\hat{\boldsymbol \xi}[t_0] \\ \hat{\m z}(k),\hat{\m p}(k)}} \;\;& \hat{\Gamma}\left(\hat{\m z}(k),\hat{\m p}(k)  \middle|\: \hat{\m x}(t_0), \left\lbrace \hat{\m d}(k) \right\rbrace_{k=t_0}^{k=t_0+H_p} \right) \notag \\
			\mathrm{s.t.}\;\;& \textit{DAE-GP}~\eqref{equ:gp-abcstracted}, \; \textit{Constraints-GP}~\eqref{equ:constr-gp-physical} \label{equ:GP-MPC} \\
			&\textit{Safety-GP}~\eqref{eq:constraints-gp-obj1},\;\textit{Smoothness-GP}~\eqref{eq:constraints-gp-obj2}. \notag
		\end{align}
	In~\eqref{equ:GP-MPC}, two sets of optimization variables are included. The first set comprises $\hat{\m x}$, $\hat{\m l}$, $\hat{\m u}$, $\hat{\m v}$, and $\hat{\m s}$ which are collected in variable $\hat{\boldsymbol \xi}[t_0]$. The latter is the corresponding GP form of $\boldsymbol \xi [t_0]$ defined in~Section~\ref{sec:CtrlObj}, that is,  $\hat{\boldsymbol \xi}[t_0] \triangleq b^{\boldsymbol \xi [t_0]}$. The second set includes the auxiliary variables $\hat{\m z}$ and $\hat{\m p}$ introduced before.
	Notice that the flow $\hat{q}_{ij}$ is an optimization variable while $q_{ij}$ is not in \eqref{equ:GP-MPC}, but a value used to calculate $C^{\mathrm{P}}(k)$, $C_{1}^{\mathrm{M}}(k)$ and $C_{2}^{\mathrm{M}}(k)$. Note that similar to the nonconvex problem~\eqref{equ:NMPC}, Problem~\eqref{equ:GP-MPC} also requires the knowledge of the state $\hat{\m x}(t_0)$ as well as a prediction of the demand pattern for an entire prediction horizon $H_p$. The detailed form of \eqref{equ:GP-MPC} is given in Table~\ref{tab:models}. We next present how we integrate the GP iterations with the MPC windows.

	
	As we mentioned in Section~\ref{sec:modelconversion}, the notation $\langle{q_{ij}(k)}\rangle_{n}$ stands for the $n$-th iteration value of $q_{ij}$ at time $k$. We use the same notation system during iterations, e.g., $\langle{\hat{\m u}(k)}\rangle_{n}$ and $\langle{\hat{\boldsymbol \xi}}\rangle_{n}$ are the $n$-th iterate value of $\hat{\m u}(k)$ and $\hat{\boldsymbol \xi}$. The algorithm of GP method for a single optimization window is presented in Algorithm~\ref{alg:seq}.
	\begin{algorithm}
		\begin{algorithmic}[1]
			\STATE \textbf{Input:} $\hat{\m x}(t_0), \{ \hat{\m d}(k) \}_{k=t_0}^{k=t_0+H_p}$
			\STATE \textbf{Output:} $\hat{\boldsymbol \xi}[t_0]$
			\STATE 	Initialize $n=0$, parameters $\langle{\hat{\boldsymbol \xi}}\rangle_0$; set $\hat{\boldsymbol \xi}_\mathrm{SaveLast} = \langle{\hat{\boldsymbol \xi}}\rangle_0$
			\WHILE 	{  $\mathrm{error} \geq \mathrm{threshold}$ \textbf{OR} $n\leq \mathrm{maxIter}$}
			\STATE  $n = n + 1$
			\STATE 	Obtain $\langle{C^{\mathrm{P}}}\rangle_n$, $\langle{C_{1}^{\mathrm{M}}}\rangle_n$ and $\langle{C_{2}^{\mathrm{M}}}\rangle_n$ by $\langle{\hat{\boldsymbol \xi}}\rangle_{n-1}$
			\STATE  Solve~\eqref{equ:GP-MPC} for $\hat{\boldsymbol \xi}_{\mathrm{Sol}}$; obtain  $\mathrm{error} = ||\hat{\boldsymbol \xi}_{\mathrm{Sol}}-\hat{\boldsymbol \xi}_{\mathrm{SaveLast}}||$
			\STATE  Assign $\langle{\hat{\boldsymbol \xi}}\rangle_{n}= \hat{\boldsymbol \xi}_{\mathrm{SaveLast}} = \hat{\boldsymbol \xi}_{\mathrm{Sol}}$
			\ENDWHILE
		\end{algorithmic}
		\caption{GP method for a single optimization window.}
		\label{alg:seq}
	\end{algorithm}
	We initialize the flow $\langle{\hat{\m u}(k)}\rangle_{0}$ and $\langle{\hat{\m v}(k)}\rangle_{0}$ in $\langle{\hat{\boldsymbol \xi}(k)}\rangle_{0}, k \in [t_0,t_0+H_p]$ with the historical statistical average flow in the pipes and pumps, and $\langle{\hat{\m s}(k)}\rangle_{0}$ is set to $\m 1$. The parameters $\langle{C^{\mathrm{P}}(k)}\rangle_1$, $\langle{C_{1}^{\mathrm{M}}(k)}\rangle_1$ and $\langle{C_{2}^{\mathrm{M}}(k)}\rangle_1$ are then calculated according to~Section~\ref{sec:modelconversion}. After solving~\eqref{equ:GP-MPC} and obtaining the $\hat{\boldsymbol \xi}_{\mathrm{Sol}}$, we save it as $\hat{\boldsymbol \xi}_\mathrm{SaveLast}$, obtain the iteration error, and make $\hat{\boldsymbol \xi}_{\mathrm{Sol}}$ as the initial value for next iteration. In addition, we define the error as the distance between current solution $\hat{\boldsymbol \xi}_{\mathrm{Sol}}$ and previous solution $\hat{\boldsymbol \xi}_{\mathrm{SaveLast}}$. The iteration continues until the error is less than a predefined error threshold ($\mathrm{threshold}$) or a maximum number of iterations ($\mathrm{maxIter}$) is reached. The details are in Algorithm~\ref{alg:seq}.
	Upon successfully solving for a single optimization window, we can continue the process and solve for the overall simulation time using a vintage MPC routine.
	
\vspace{-0.35cm}

	\section{Numerical Tests}~\label{sec:test}
	In this section, we present a simple simulation example to illustrate the applicability of the GP-based MPC formulation for WDN operations. The considered water network is an 8-node network from EPANET Users Manual~\cite[Chapter 2]{rossman2000epanet}. 
	The numerical tests are simulated using EPANET Matlab Toolkit~\cite{eliades2009epanet} with GP solver \textit{GGPLAB}~\cite{mutapcic2006ggplab}.
	The basic parameters in 8-node network including the elevation of nodes, length, and diameter of pipes are obtained from~\cite{rossman2000epanet}.  We now present the list of constraints and parameters used in the simulations. 
		\begin{figure}[h]
		\centering
		\includegraphics[width=1\linewidth]{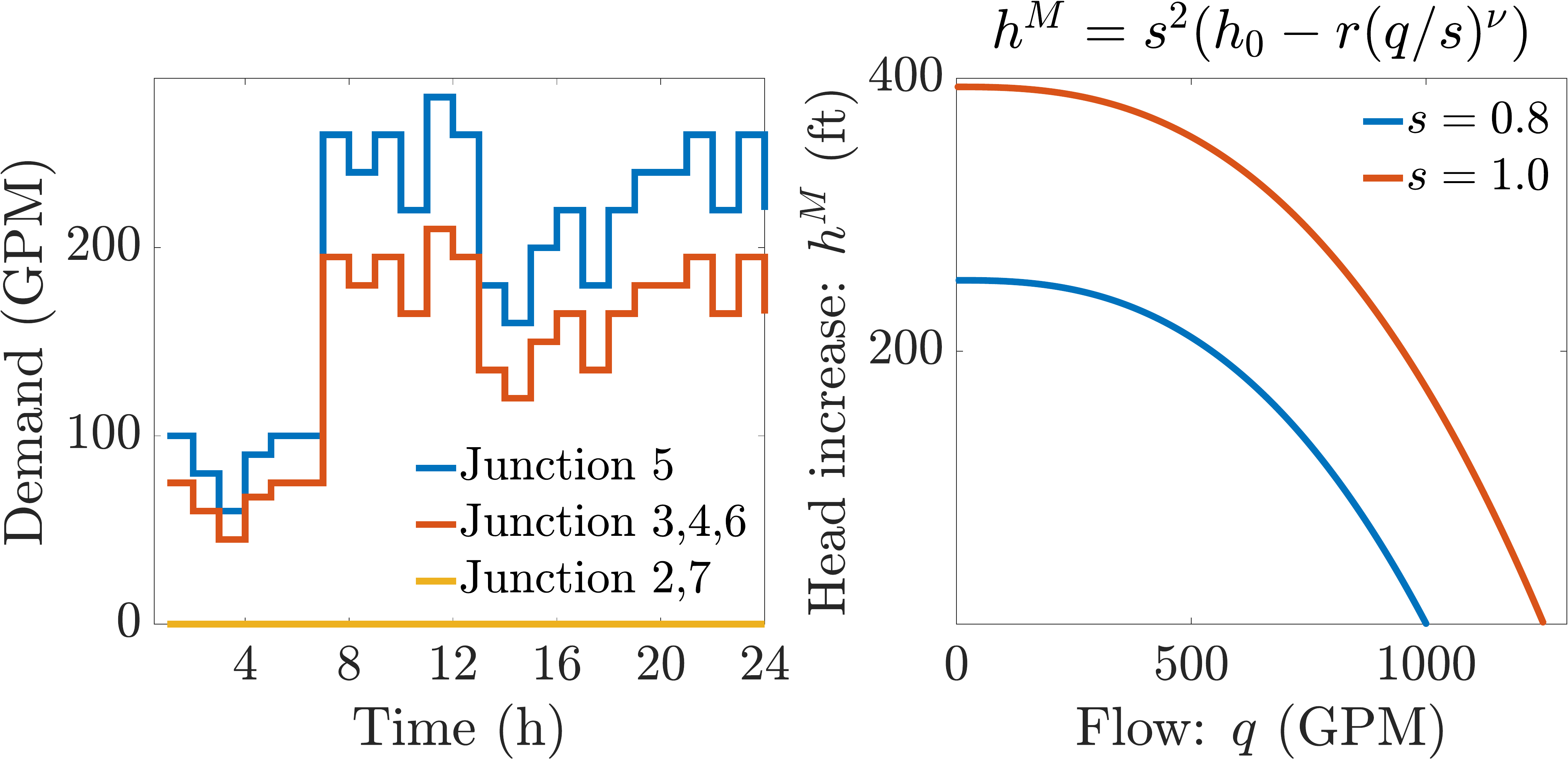}
		\caption{Water demand at various junctions and variable-speed pump curve ($h_0 = 393.7$, $r = 3.7\times 10^{-6}$,  $\nu = 2.59$).}
		\label{fig:DemandAndCurve}
		\vspace{-0.2cm}
	\end{figure}	

	\noindent $\bullet$  The initial head of Tank 8 is $834\;\mathrm{ft}$, the water level range of Tank 8 is $[830,850]\;\mathrm{ft}$, and the safety water level $x^\mathrm{sf}$~\eqref{equ:SafetyWater} of Tank 8 from Section~\ref{sec:CtrlObj} is set to $838\; \mathrm{ft}$.   We set the total simulation time $T_\mathrm{final}$ to $24\,\mathrm{hrs}$ in Algorithm~\ref{alg:seq}. The demand pattern for 24 hours at different junctions is shown in Fig.~\ref{fig:DemandAndCurve}. This demand pattern is different from~\cite{rossman2000epanet}, as our intention is to make the demand vary more rapidly to test the performance of the presented GP-based MPC routine.
	
	\noindent $\bullet$ The relationship between head increase and flow of Pump 9 defined by~\eqref{equ:head-flow-pump} is presented in~Fig.~\ref{fig:DemandAndCurve}. We observe that the head increase and flow provided by a pump varies with the relative speed $s \in [0,1]$ with $s = 0$ refers to the pump being off. In~\eqref{equ:constraints}, the physical constraints of the head imposed at the $i$-th junction is greater than its corresponding elevation, and the head of $i$-th reservoir is fixed at its elevation. Since we have only one reservoir, this implies that $h_1^{\mathrm{R}} = 700\;\mathrm{ft}$. As for the flow, the direction is unknown, and we simply constrain the flow to $q_i \in [-1000,1000]\;\mathrm{GPM}$.
	
	\noindent $\bullet$  For the geometric programming component of the presented formulations, we set the base $b = 1.005$. The parameters we use in Algorithm~\ref{alg:seq} are selected as: $\mathrm{error} = 0.5$ and $\mathrm{maxIter}= 40$. We consider a sampling time of $1\,\mathrm{hr}$, a prediction horizon $H_p = 10\,\mathrm{hrs}$, and a simulation of $24\,\mathrm{hrs.}$  

\begin{figure}[h]
		\centering
		\includegraphics[width=1\linewidth]{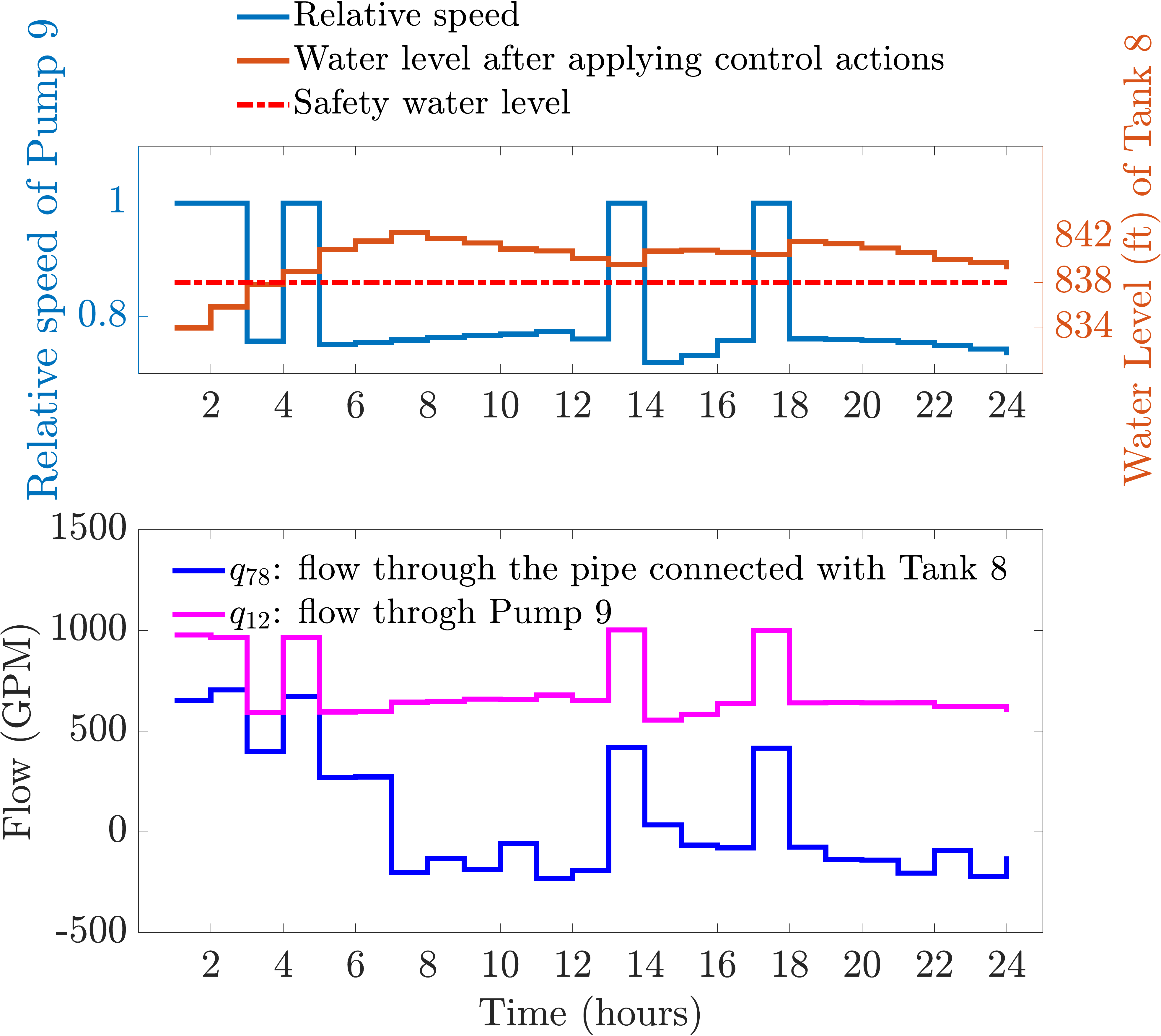}
		\vspace{-0.4cm}
		\caption{Relative speed of Pump 9, controlled water level of Tank 8, and flows through a pipe and Pump 9.}
		\label{fig:WaterLever_PumpSpeed_Flow}
	\end{figure}

	Here, we present the solution to WDNs operation problem after running Algorithm~\ref{alg:seq}. Fig.~\ref{fig:WaterLever_PumpSpeed_Flow} shows the optimal control effort (the variable pump speed), the water level after applying the pump control, and the water flow through some pipes. During time period $k \in [1,2]$, the water level of Tank 8 is below the safety water level $838\;\mathrm{ft}$, so the relative speed of Pump 9 is set to the highest speed $s=1$ by GP-MPC controller. Consequently,  Pump 9 starts to deliver water into the network and more water begins to flow into the tank. At $k=3$, the target setpoint is reached and the pump slows down to a relative speed of $s=0.7566$ according to MPC algorithm. Notice that the demand at $k=7$ increases dramatically as shown in Fig.~\ref{fig:DemandAndCurve}, and as a result the pump speed increases to the maximum speed at $k=4$ to prepare for this situation in advance. 
	\begin{figure}[h]
		\vspace{-0.45cm}
		\centering
		\includegraphics[width=1\linewidth]{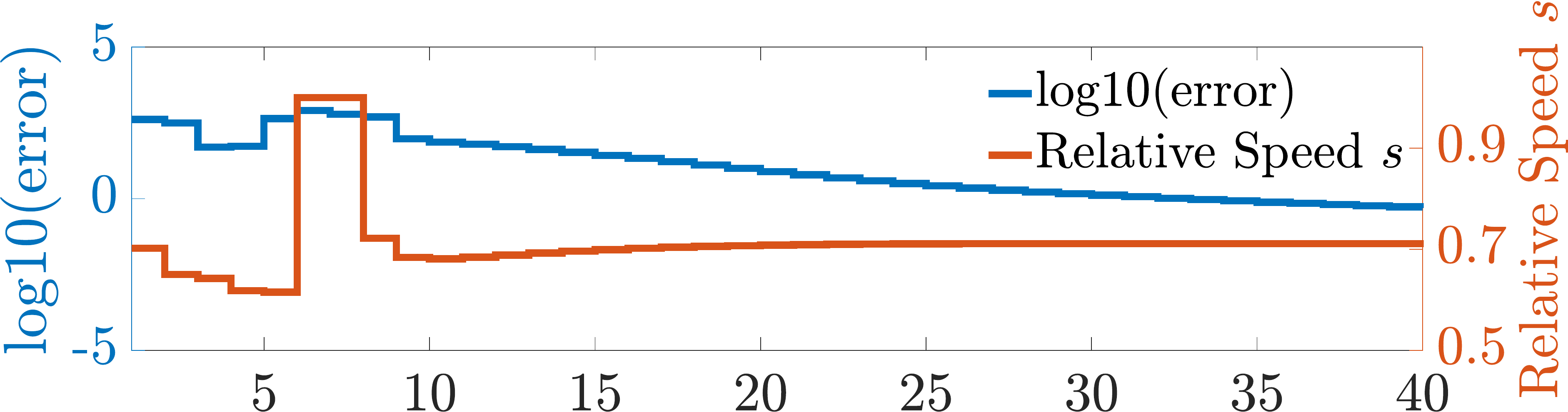}
		\vspace{-0.34cm}
		\caption{Convergence of the Iteration error and relative speed from Algorithm~\ref{alg:seq}.}
		\label{fig:Convergence}
	\end{figure}
	During time period $k \in [7,13]$, the demand changes slightly, and the relative speed is stabilized around $s=0.8$. At $k=14$, the water level decreases to $838\;\mathrm{ft}$, which leads to an increase in the pump speed to maintain the desired, safe water level. During time period $k \in [18,24]$, the increase in demand justifies the higher pump speed which is then followed by a decrease in the speed for the last stretch of the day. 
	Fig.~\ref{fig:Convergence} shows the convergence of Algorithm~\ref{alg:seq}---the GP method for a single optimization window when $t_0=23$. We can see the iteration error keeps decreasing which implies the  final solution is reached slowly. For the first 7 iterations, the ups and downs of error is caused by whether the objective function~\eqref{eq:gp-obj} is activated, because the solved $\m x$ from Algorithm~\ref{alg:seq} is not accurate and varies. When the solved $\m x$ becomes stable, the objective function~\eqref{eq:gp-obj} does not change, and Algorithm~\ref{alg:seq} converges to the final solution. 
	It is noteworthy to mention that the flow $q_{12}(k)$ through Pump 9 is always in one direction, while the flow direction of $q_{78}(k)$ through the pipe connected to Tank 8 changes 5 times during the simulation in Fig.~\ref{fig:WaterLever_PumpSpeed_Flow}. This justifies the use of the \textit{flow direction-unaware GP} component in the MPC. Future work will experiment with larger networks and investigate the impact of demand uncertainty on the WDN control.

	\bibliographystyle{IEEEtran}
	\bibliography{IEEEabrv,bibfile2}

\begin{thebibliography}{10}
\providecommand{\url}[1]{#1}
\csname url@samestyle\endcsname
\providecommand{\newblock}{\relax}
\providecommand{\bibinfo}[2]{#2}
\providecommand{\BIBentrySTDinterwordspacing}{\spaceskip=0pt\relax}
\providecommand{\BIBentryALTinterwordstretchfactor}{4}
\providecommand{\BIBentryALTinterwordspacing}{\spaceskip=\fontdimen2\font plus
\BIBentryALTinterwordstretchfactor\fontdimen3\font minus
  \fontdimen4\font\relax}
\providecommand{\BIBforeignlanguage}[2]{{%
\expandafter\ifx\csname l@#1\endcsname\relax
\typeout{** WARNING: IEEEtran.bst: No hyphenation pattern has been}%
\typeout{** loaded for the language `#1'. Using the pattern for}%
\typeout{** the default language instead.}%
\else
\language=\csname l@#1\endcsname
\fi
#2}}
\providecommand{\BIBdecl}{\relax}
\BIBdecl

\bibitem{cross1936analysis}
H.~Cross, ``Analysis of flow in networks of conduits or conductors,''
  University of Illinois at Urbana Champaign, College of Engineering.
  Engineering Experiment Station., Tech. Rep., 1936.

\bibitem{wood1972hydraulic}
D.~J. Wood and C.~O. Charles, ``Hydraulic network analysis using linear
  theory,'' \emph{Journal of the Hydraulics division}, vol.~98, no.~7, pp.
  1157--1170, 1972.

\bibitem{arora1976flows}
M.~L. Arora, ``Flows split in closed loops expending least energy,''
  \emph{Journal of the Hydraulics Division}, vol. 102, no.~3, pp. 455--458,
  1976.

\bibitem{Hafez-FixedPointWDSA}
\BIBentryALTinterwordspacing
M.~Bazrafshan, N.~Gatsis, M.~Giacomoni, and A.~Taha, ``A fixed-point iteration
  for steady-state analysis of water distribution networks,'' in \emph{Proc.
  6th IEEE Global Conf. Signal and Information Processing}, Anaheim, CA, Nov.
  2018.  \url{https://arxiv.org/abs/1807.01404}
\BIBentrySTDinterwordspacing

\bibitem{singh2019flow}
M.~K. Singh and V.~Kekatos, ``On the flow problem in water distribution
  networks: Uniqueness and solvers,'' \emph{arXiv preprint arXiv:1901.03676},
  2019.

\bibitem{mala2017lost}
H.~Mala-Jetmarova, N.~Sultanova, and D.~Savic, ``Lost in optimisation of water
  distribution systems? a literature review of system operation,''
  \emph{Environmental Modelling \& Software}, vol.~93, pp. 209--254, 2017.

\bibitem{wang2017non}
Y.~Wang, V.~Puig, and G.~Cembrano, ``Non-linear economic model predictive
  control of water distribution networks,'' \emph{Journal of Process Control},
  vol.~56, pp. 23--34, 2017.

\bibitem{sankar2015optimal}
G.~S. Sankar, S.~M. Kumar, S.~Narasimhan, S.~Narasimhan, and S.~M. Bhallamudi,
  ``Optimal control of water distribution networks with storage facilities,''
  \emph{Journal of Process Control}, vol.~32, pp. 127--137, 2015.

\bibitem{wang2016stochastic}
Y.~Wang, C.~Ocampo-Martinez, and V.~Puig, ``Stochastic model predictive control
  based on gaussian processes applied to drinking water networks,'' \emph{IET
  Control Theory \& Applications}, vol.~10, no.~8, pp. 947--955, 2016.

\bibitem{sun2016combining}
C.~C. Sun, V.~Puig, and G.~Cembrano, ``Combining csp and mpc for the
  operational control of water networks,'' \emph{Engineering Applications of
  Artificial Intelligence}, vol.~49, pp. 126--140, 2016.

\bibitem{humpola2013unified}
J.~Humpola and A.~F{\"u}genschuh, ``A unified view on relaxations for a
  nonlinear network flow problem,'' 2013.

\bibitem{singh2018optimal}
M.~K. Singh and V.~Kekatos, ``Optimal scheduling of water distribution
  systems,'' \emph{arXiv preprint arXiv:1806.07988}, 2018.

\bibitem{xie2015nonlinear}
M.~Xie and M.~Brdys, ``Nonlinear model predictive control of water quality in
  drinking water distribution systems with dbps objectives.''

\bibitem{sela2015control}
L.~Sela~Perelman and S.~Amin, ``Control of tree water networks: A geometric
  programming approach,'' \emph{Water Resources Research}, vol.~51, no.~10, pp.
  8409--8430, 2015.

\bibitem{fooladivanda2017energy}
D.~Fooladivanda and J.~A. Taylor, ``Energy-optimal pump scheduling and water
  flow,'' \emph{IEEE Transactions on Control of Network Systems}, no.~3, pp.
  1016--1026, Sept. 2018.

\bibitem{zamzam2018optimal}
A.~S. Zamzam, E.~Dall'Anese, C.~Zhao, J.~A. Taylor, and N.~Sidiropoulos,
  ``Optimal water-power flow problem: Formulation and distributed optimal
  solution,'' \emph{IEEE Transactions on Control of Network Systems}, 2018.

\bibitem{rossman2000epanet}
L.~A. Rossman \emph{et~al.}, ``Epanet 2: users manual,'' 2000.

\bibitem{ocampo2013application}
C.~Ocampo-Martinez, V.~Puig, G.~Cembrano, and J.~Quevedo, ``Application of
  predictive control strategies to the management of complex networks in the
  urban water cycle applications of control.''\hskip 1em plus 0.5em minus
  0.4em\relax Institute of Electrical and Electronics Engineers, 2013.

\bibitem{boyd2007tutorial}
S.~Boyd, S.-J. Kim, L.~Vandenberghe, and A.~Hassibi, ``A tutorial on geometric
  programming,'' \emph{Optimization and engineering}, vol.~8, no.~1, p.~67,
  2007.

\bibitem{eliades2009epanet}
D.~El{\'\i}ades and M.~Kyriakou, ``Epanet matlab toolkit,'' \emph{University of
  Cyprus, Republic of Cyprus}, 2009.

\bibitem{mutapcic2006ggplab}
A.~Mutapcic, K.~Koh, S.~Kim, L.~Vandenberghe, and S.~Boyd, ``Ggplab: A simple
  matlab toolbox for geometric programming,'' \emph{web page and software:
  http://stanford. edu/boyd/ggplab}, 2006.

\end{thebibliography}

	%
\end{document}